\documentclass{amsart}
\usepackage{graphicx}
\calclayout
\renewcommand\section[1]{\medbreak\par{\bf#1}\medskip\nobreak}

\begin{document}
\title{Nefarious Numbers}
\author{Douglas N.~Arnold and Kristine K.~Fowler}
\thanks{Douglas N.~Arnold is McKnight Presidential Professor of Mathematics at the University of Minnesota and president of the
Society for Industrial and Applied Mathematics. Kristine K.~Fowler is mathematics librarian at the University of Minnesota. The
authors gratefully acknowledge the assistance of Susan K.~Lowry, who developed and supported the database used in this study,
and Molly T.~White.}
\maketitle
\thispagestyle{empty}

\section{Introduction}

The impact factor has been widely adopted as a proxy for journal quality.  It is used by libraries to guide purchase and renewal decisions, by researchers deciding where to publish and what to read, by tenure and promotion committees laboring under the assumption that publication in a higher impact factor journal represents better work, and by editors and publishers as a means to evaluate and promote their journals. The impact factor for a journal in a given year is calculated by ISI (Thomson Reuters) as the average number of citations in that year to the articles the journal published in the preceding two years.  It has been widely criticized on a variety of grounds\footnote{P.~O.~Seglen, Why the impact factor of journals should not be used for evaluating research. BMJ 314 (1997), 498--502.}${}^,$%
\footnote{J.~Ewing, Measuring journals. Notices of the AMS 53 (2006), 1049--1053.}${}^,$%
\footnote{R.~Golubic, M.~Rudes, N.~Kovacic, M.~Marusic, and A.~Marusic, Calculating impact factor: how bibliographical classification of
journal items affects the impact factor of large and small journals. Sci. Eng. Ethics 14 (2008), 41--49.}${}^,$%
\footnote{R.~Adler, J.~Ewing, and P.~Taylor, Citation statistics. Statistical Sciences 24 (2009), 1--14.}:

\begin{itemize}

\item A journal's distribution of citations does not determine its quality.

\item The impact factor is a crude statistic, reporting only one particular item of information from the citation distribution.

\item It is a flawed statistic.  For one thing, the distribution of citations among papers is highly skewed, so the mean for the journal tends to be misleading.  For another, the impact factor only refers to citations within the first two years after publication (a particularly serious deficiency for mathematics, in which around 90\% of citations occur after two years).

\item The underlying database is flawed, containing errors and including a biased selection of journals.

\item Many confounding factors are ignored, for example, article type (editorials, reviews, and letters versus original research articles), multiple authorship, self-citation, language of publication, etc.

\end{itemize}

Despite these difficulties, the allure of the impact factor as a single, readily available number---not requiring complex judgments or expert input, but purporting to represent journal quality---has proven irresistible to many.  Writing in 2000 in a newsletter for journal editors, Amin and Mabe\footnote{M.~Amin and M.~Mabe, Impact factors: use and abuse. Perspectives in Publishing 1 (2000), 1--6.} wrote that the ``impact factor has moved in recent years from an obscure bibliometric indicator to become the chief quantitative measure of the quality of a journal, its research papers, the researchers who wrote those papers and even the institution they work in.'' It has become commonplace for journals to issue absurd announcements touting their impact factors, like this one which was mailed around the world by World Scientific, the publisher of the International Journal of Algebra and Computation: ``IJAC's Impact Factor has improved from 0.414 in 2007 to 0.421 in 2008!  Congratulations to the Editorial Board and contributors of IJAC.'' In this case, the 1.7\% increase in the impact factor represents a single additional citation to one of the 145 articles published by the journal in the preceding two years.

Because of the (misplaced) emphasis on impact factors, this measure has become a target at which journal editors and publishers aim.  This has in turn led to another major source of problems with the factor.  Goodhart's law warns us that ``when a measure becomes a target, it ceases to be a good measure.''\footnote{This succinct formulation is from M.~Strathern, `Improving ratings': audit in the British University system, European Review 5
(1997), 305--321.}  This is precisely the case for impact factors.  Their limited utility has been further compromised by impact factor manipulation, the engineering of this supposed measure of journal quality, in ways that increase the measure, but do not add to---indeed subtract from---journal quality.

Impact factor manipulation can take numerous forms.  In a 2007 essay on the deleterious effects of impact factor manipulation, Macdonald and Kam\footnote{S.~Macdonald and J.~Kam, Aardvark et al.: quality journals and gamesmanship in management studies. Journal of Information Science 33 (2007), 702--717.} noted wryly that ``the canny editor cultivates a cadre of regulars who can be relied upon to boost the measured quality of the journal by citing themselves and each other shamelessly.'' There have also been widespread complaints by authors of manuscripts under review, who were asked or required by editors to cite other papers from the journal.  Given the dependence of the author on the editor's decision for publication, this practice borders on extortion, even when posed as a suggestion.  In most cases, one can only guess about the presence of such pressures, but overt instances were reported already in 2005 by Monastersky\footnote{R.~Monastersky, The number that's devouring science. Chronicle of Higher Education 52 (2005).} in the Chronicle of Higher Education and Begley\footnote{S.~Begley, Science journals artfully try to boost their rankings. Wall Street Journal, 5 June 2006, B1.} in the Wall Street Journal. A third well-established technique by which editors raise their journals' impact factors, is by publishing review items with large numbers of citations to the journal.  For example, the Editor-in-Chief of the Journal of Gerontology A made a practice of authoring and publishing a review article every January focusing on the preceding two years; in 2004, 195 of the 277 references were to the Journal of Gerontology A. Though the distortions these unscientific practices wreak upon the scientific literature have raised occasional alarms, many suppose that they either have minimal effect or are so easily detectable they can be disregarded. A counterexample should confirm the need for alarm.

\section{The case of IJNSNS}

The field of applied mathematics provides an illuminating case in which we can study such impact factor distortion.  For the last several years, the International Journal of Nonlinear Sciences and Numerical Simulation (IJNSNS) has dominated the impact factor charts in the ``Mathematics, Applied'' category.  It took first place in each year 2006, 2007, 2008, and 2009, generally by a wide margin, and came in second in 2005.  However, as we shall see, a more careful look indicates that IJNSNS is nowhere near the top of its field.  Thus we set out to understand the origin of its large impact factor.

In 2008, the year we shall consider in most detail, IJNSNS had an impact factor of 8.91, easily the highest among the 175 journals in the applied math category in ISI's Journal Citation Reports (JCR).  As controls, we will also look at the two journals in the category with the second and third highest impact factors, Communications on Pure and Applied Mathematics (CPAM), and SIAM Review (SIREV), with 2008 impact factors of 3.69 and 2.80, respectively.  CPAM is closely associated with the Courant Institute of Mathematical Sciences, and SIREV is the flagship journal of the Society for Industrial and Applied Mathematics (SIAM).\footnote{The first author is the current president of SIAM.}  Both journals have a reputation for excellence.
Evaluation based on expert judgment is the best alternative to citation-based measures for journals.  Though not without potential problems of its own, a careful rating by experts is likely to provide a much more accurate and holistic guide to journal quality than impact factor or similar metrics.  In mathematics, as in many fields, researchers are widely in agreement about which are the best journals in their specialties.  The Australian Research Council recently released such an evaluation, listing quality ratings for over 20,000 peer-reviewed journals across disciplines. The list was developed through an extensive review process involving learned academies (such as the Australian Academy of Science), disciplinary bodies (such as the Australian Mathematical Society), and many researchers and expert reviewers.\footnote{Australian Research Council, Ranked Journal List Development, http://www.arc.gov.au/era/journal\_list\_dev.htm.}  This rating will be used in 2010 for the Excellence in Research Australia assessment initiative, and is referred to as the ERA 2010 Journal List.  The assigned quality rating, which is intended to represent ``the overall quality of the journal,'' is one of four values:

\begin{itemize}

\item A*: one of the best in its field or subfield
\item A: very high quality
\item B: solid, though not outstanding reputation
\item C: does not meet the criteria of the higher tiers

\end{itemize}

The ERA list included all but five of the 175 journals assigned a 2008 impact factor by JCR in the category
\begin{figure}[ht]
\centerline{\includegraphics[height=5in]{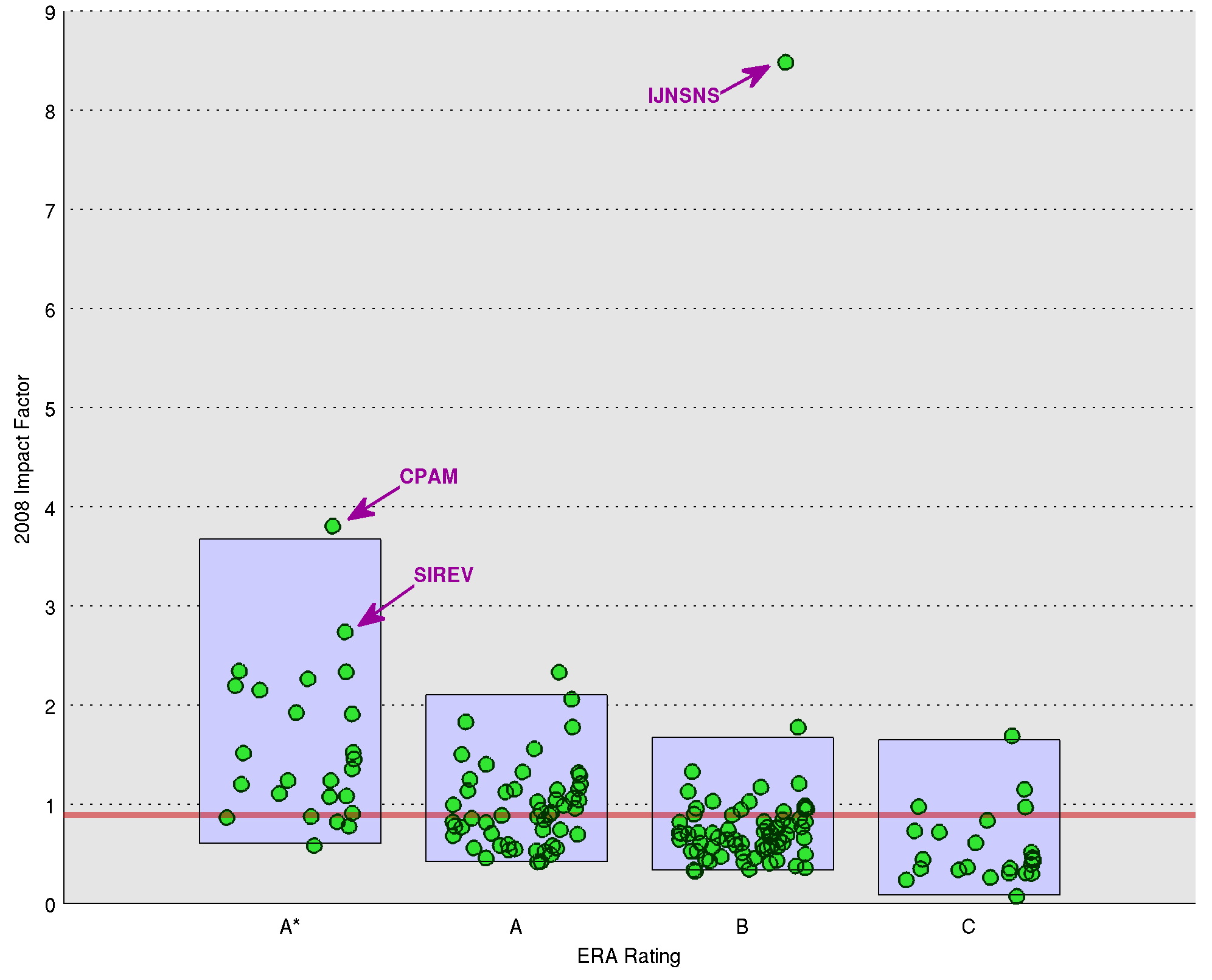}}
\caption{\it 2008 impact factors of 170 applied math journals grouped according to their 2010 ERA rating tier. In each tier, the band runs from the 2.5th to the 97.5th percentile, outlining the middle 95\%. Horizontal position of the data points within tiers is assigned randomly to improve visibility. The red line is at the 20th percentile of the $A^*$ tier.}
\end{figure}
``Mathematics, Applied.''  Figure 1 shows the impact factors for journals in each of the four rating tiers.
We see that, as a proxy for expert opinion, the impact factor does rather poorly.  There are many examples of journals with a higher impact factor than other journals which are one, two, and even three rating tiers higher.  The red line is drawn
so that 20\% of the A* journals are below it; it is notable that 51\% of the A journals have an impact factor above that level, as do 23\% of the B journals and even 17\% of those in the C category. The most extreme outlier is IJNSNS, which, despite its relatively astronomical impact factor, is not in the first or second, but rather third tier.  The ERA rating assigned its highest score, A*, to 25 journals.  Most of the journals with the highest impact factors are here, including CPAM and SIREV, but of the top 10 journals by impact factor, two were assigned an A, and only IJNSNS was assigned a B.  There were 53 A-rated journals, and 69 B-rated journals altogether.  If IJNSNS were assumed to be the best of the B journals, there would be 78 journals with higher ERA ratings, while if it were the worst, its ranking would fall to 147.  In short, the ERA ratings suggest that IJNSNS is not only not the top applied math journal, but its rank should be somewhere in the range 75--150.  This remarkable mismatch between reputation and impact factor begs an explanation.

\section{Makings of a high impact factor}

A first step to understanding IJNSNS's high impact factor is to look at how many authors contributed substantially to the counted citations, and who they were.  The top-citing author to IJNSNS in 2008 was the journal's Editor-in-Chief, Ji-Huan He, who cited the journal (within the two-year window) 243 times.  The second top-citer, D.~D.~Ganji, with 114 cites, is also a member of the editorial board, as is the third, regional editor Mohamed El Naschie, with 58 cites.  Together these three account for 29\% of the citations counted towards the impact factor.  For comparison, the top three citers to SIREV contributed only 7, 4, and 4 citations, respectively, accounting for less than 12\% of the counted citations, and none of these authors is involved in editing the journal.  For CPAM the top three citers (9, 8, and 8) contributed about 7\% of the citations, and, again, were not on the editorial board.
Another significant phenomenon is the extent to which citations to IJNSNS are concentrated within the 2-year window used in the impact factor calculation.  Our analysis of 2008 citations to articles published since 2000 shows that 16\% of the citations to CPAM fell within that 2-year window, and only 8\% of those to SIREV did; in contrast, 71.5\% of the 2008 citations to IJNSNS fell within the 2-year window.  In Table 1, we show the 2008 impact factors for the three journals, as well as a modified impact factor, which gives the average number of citations in 2008 to articles the journals published not in 2006 and 2007, but in the preceding six years.  Since the cited half-life (the time it takes to generate half of all the eventual citations to an article) for applied mathematics is nearly 10 years,\footnote{In 2010, Journal Citation Reports assigned the category ``Mathematics, Applied'' an aggregate cited half-life of 9.5 years.} this measure is at least as reasonable as the impact factor. It is also independent, unlike JCR's 5-Year Impact Factor, as its time period does not overlap with that targeted by the impact factor.
\begin{table}[ht]
 \centerline{\begin{tabular}{|c|c|c|}
\hline
&
2008 impact factor with&
Modified 2008 ``impact factor''\\
Journal& normal 2006--7 window& with 2000--5 window\\
\hline
IJNSNS&
8.91&
1.27\\
\hline
CPAM&
3.69&
3.46\\
\hline
SIREV&
2.8&
10.4\\
\hline
\end{tabular}}
\vspace{10pt}
\caption{\it 2008 impact factors computed with the usual two-preceding years window, and with a window going back eight years but neglecting the two immediately preceding.}
\end{table}
Note that the impact factor of IJNSNS drops precipitously, by a factor of seven, when we consider a different citation window.  By contrast the impact factor of CPAM stays about the same and that of SIREV increases markedly.  One may simply note that, in distinction to the controls, the citations made to IJNSNS in 2008 greatly favor articles published in precisely the two years which are used to calculate the impact factor.

Further striking insights arise when we examine the high-citing journals rather than high-citing authors.  The counting of journal self-citations in the impact factor is frequently criticized, and indeed it does come into play in this case.  In 2008, IJNSNS supplied 102, or 7\%, of its own impact factor citations.  The corresponding numbers are 1 citation (0.8\%) for SIREV and 8 citations (2.4\%) for CPAM.  The disparity in other recent years is similarly large or larger.

However, it was Journal of Physics: Conference Series, which provided the greatest number of IJNSNS citations.  A single issue of that journal provided 294 citations to IJNSNS in the impact-factor window, accounting for more than 20\% of its impact factor.  What was this issue?  It was the proceedings of a conference organized by IJNSNS Editor-in-Chief He at his home university. He was responsible for the peer review of the issue. The second top-citing journal for IJNSNS was Topological Methods in Nonlinear Analysis, which contributed 206 citations (14\%), again with all citations coming from a single issue.  This was a special issue with Ji-Huan He as the guest editor; his co-editor, Lan Xu, is also on the IJNSNS editorial board.  J.-H.~He himself contributed a brief article to the special issue, consisting of 3 pages of text and 30 references.  Of these, 20 were citations to IJNSNS within the impact-factor window.  The remaining 10 consisted of 8 citations to He and 2 to Xu.

Continuing down the list of IJNSNS high-citing journals, another similar circumstance comes to light: 50 citations from a single issue of the Journal of Polymer Engineering (which, like IJNSNS, is published by Freund), guest-edited by the same pair Ji-Huan He and Lan Xu.  However, third place is held by the journal Chaos, Solitons \& Fractals, with 154 citations spread over numerous issues.  These are again citations which may be viewed as subject to editorial influence or control.  In 2008 Ji-Huan He served on the editorial board of CS\&F, and its Editor-in-Chief was Mohamed El Naschie, who was also a co-editor of IJNSNS. In a highly publicized case, the entire editorial board of CS\&F was recently replaced, but El Naschie remained co-editor of IJNSNS.

Many other citations to IJNSNS came from papers published in journals for which He served as editor, such as Zeitschrift f\"ur Naturforschung A, which provided 40 citations; there are too many others to list here, since He serves in an editorial capacity on more than 20 journals (and has just been named Editor-in-Chief of four more journals from the newly-formed Asian Academic Publishers). Yet another source of citations came from papers authored by IJNSNS editors other than He, which accounted for many more.  All told, the aggregation of such editor-connected citations, which are time-consuming to detect, account for more than 70\% of all the citations contributing to the IJNSNS impact factor.

\section{Bibliometrics for individuals}

Bibliometrics are also used to evaluate individuals, articles, institutions and even nations.  Essential Science Indicators, which is produced by Thomson Reuters, is promoted as a tool for ranking ``top countries, journals, scientists, papers, and institutions by field of research.''  However, these metrics are primarily based on the same citation data used for journal impact factors and thus they can be manipulated just as easily, indeed simultaneously.  The special issue of Journal of Physics: Conference Series which He edited and which garnered 243 citations for his journal, also garnered 353 citations to He himself. He claims a total citation count of over 6,800.\footnote{This claim, and that of an h-index of 39, are made in the biographical notes of one of his recent papers (Nonl.~Sci.~Letters~1
(2010), page 1).} Even half that is considered highly noteworthy as evidenced by this announcement in ScienceWatch.com:\footnote{ScienceWatch.com, April 2008, http://sciencewatch.com/inter/aut/2008/08-apr/08aprHe/.} ``According to a recent analysis of Essential Science Indicators from Thomson Scientific, Professor Ji-Huan He has been named a Rising Star in the field of Computer Science... His citation record in the Web of Science includes 137 papers cited a total of 3,193 times to date.''  Together with only a dozen other scientists in all fields of science, He was cited by ESI for the ``Hottest Research of 2007--8'' and again for the ``Hottest Research of 2009.''

The h-index is another popular citation-based metric for researchers, intended to measure productivity as well as impact.  An individual's h-index is the largest number such that that many of his or her papers have been cited at least that many times. It too is not immune from Goodhart's law.  J.-H.~He claims an h-index of 39, while Hirsch estimated the median for Nobel prize winners in physics to be 35.\footnote{J.~Hirsch, An index to quantify an individual's scientific research output. PNAS 102 (2005), 16569--16572.}
Whether for judgment of individuals or journals, citation-based designations are no substitute for an informed judgment of quality.

\newpage

\section{Closing thoughts}

Despite numerous flaws, the impact factor has been widely used as a measure of quality for journals, and even for papers and authors.  This creates an incentive to manipulate it.  Moreover, it is possible to vastly increase impact factor without increasing journal quality at all. The actions of a few interested individuals can make a huge difference, yet require considerable digging to reveal.  We primarily discussed one extreme example, but there is little reason to doubt that such techniques are being used to a lesser---and therefore less easily detected---degree by many journals.  The cumulative result of the design flaws and manipulation is that impact factor gives a very inaccurate view of journal quality.  More generally, the citations which form the basis of the impact factor and various other bibliometrics are inherently untrustworthy.

The consequences of this unfortunate situation are great.  Rewards are wrongly distributed, the scientific literature and enterprise are distorted, and cynicism about them grows.  What is to be done?  Just as for scientific research itself, the temptation to embrace simplicity when it seriously compromises accuracy, must be resisted.  Scientists who give in to the temptation to suppress data or fiddle with statistics to draw a clearer point are censured.  We must bring a similar level of integrity to the evaluation of research products.  Administrators, funding agencies, librarians, and others needing such evaluations should just say no to simplistic solutions, and approach important decisions with thoughtfulness, wisdom, and expertise.

\end{document}